\documentclass[12pt]{article}
\usepackage{amsmath,amssymb,amsbsy,amsfonts,amsthm,
            latexsym,amsopn,amstext,amsxtra,euscript,amscd}

\begin{document}

\newtheorem{prop}{Proposition}
\newtheorem{lemma}[prop]{Lemma}
\newtheorem{cor}[prop]{Corollary}
\newtheorem{theorem}[prop]{Theorem}
\newtheorem{com}[prop]{Comments}
\newtheorem{remark}[prop]{Remark}

\title{\bf On uniform distribution modulo one}
\author{
{\sc M.~Z.~Garaev}\\
{\normalsize Instituto\ de Matem{\'a}ticas,
UNAM} \\
{\normalsize Campus Morelia, Ap. Postal 61-3 (Xangari)}\\
{\normalsize C.P.\ 58089, Morelia, Michoac{\'a}n, M{\'e}xico} \\
{\normalsize \tt garaev@matmor.unam.mx}}

\date{}

\maketitle

\begin{abstract}
We introduce an elementary argument to the theory of distribution of
sequences modulo one.
\end{abstract}


\section{Introduction}

Throughout the paper $x_1,x_2,\ldots $ denotes  a sequence of real
numbers with their fractional parts $\{x_1\}, \{x_2\},\ldots.$ For
$0\le\alpha<\beta\le 1$ we use $F(N, x_n; \alpha,\beta)$ to denote
the number of terms of this sequence with the condition
$$
\alpha\le \{x_n\}<\beta, \quad n\le N.
$$
The sequence $x_n$ is called uniformly distributed modulo one if
$$
\lim_{N\to\infty}\ \sup_{0\le \alpha<\beta\le 1}\left|\frac{F(N,
x_n; \alpha,\beta)}{N}-(\beta-\alpha)\right|=0.
$$
The central place in the theory of uniform distribution modulo one
belongs to the Weyl criterion. Its most nontrivial part reads as
follows: {\it if for any integer $h\not=0$ we have
$$
\lim_{N\to\infty}\frac{1}{N}\sum_{n=1}^Ne^{2\pi ihx_n}=0,
$$
then $x_n$ is uniformly distributed modulo one.} It is easy to see
that the opposite statement is also true.

The traditional method to obtain quantified versions of the Weyl
criterion is Vinogradov's lemma on ``little glasses'', see
Vinogradov~\cite[Lemma 2, Chapter II]{Vi} or Karatsuba~\cite[Lemma
A, Chapter I]{Ka}. The well known Erd\H os-Tur{\'a}n inequality
claims that for any $H\ge 1,$
$$
\sup_{0\le \alpha<\beta\le 1}\left|\frac{F(N, x_n;
\alpha,\beta)}{N}-(\beta-\alpha)\right| \ll
\frac{1}{H}+\frac{1}{N}\sum_{h=1}^{H}\frac{1}{h}\left|\sum_{n=1}^Ne^{2\pi
ihx_n}\right|,
$$
see Montgomery~\cite[Corollary~1.1, Chapter I]{Monty}.
In~\cite[Theorem 1, Chapter 1]{Monty} the following estimate has
been proved:
\begin{equation}
\label{eqn:Monty} \left|\frac{F(N, x_n;
\alpha,\beta)}{N}-(\beta-\alpha)\right| \ll
\frac{1}{H}+\frac{1}{N}\sum_{h=1}^{H}\min\left(\beta-\alpha,
\frac{1}{h}\right)\left|\sum_{n=1}^Ne^{2\pi ihx_n}\right|.
\end{equation}
The advantage of~\eqref{eqn:Monty} over the Erd\H os-Tur{\'a}n
inequality is that it gives more precise information on distribution
of $\{x_n\}$ in small intervals.

The aim of the present paper is to introduce an elementary
self-contained argument to investigate the problem of uniform
distribution of sequences modulo one.

Throughout the paper we use the following simple identity:
\begin{equation*}
\label{eq:Ident} \frac{1}{m}\sum_{h=0}^{m-1} e^{2\pi i h\frac{u}{m}}
= \left\{
\begin{array}{ll}
0,& \quad \mbox{if}\ u\not \equiv 0 \pmod m, \\
1,& \quad \mbox{if}\ u \equiv 0 \pmod m.
\end{array} \right.
\end{equation*}
In particular, if $\mathcal{X}\in\{0,1,\ldots,m-1\}$ is a set with
$|\mathcal{X}|$ elements, then
$$
\frac{1}{m}\sum_{h=0}^{m-1} \left|\sum_{x\in\mathcal{X}}e^{2\pi i
h\frac{x}{m}}\right|^2=\frac{1}{m}\sum_{h=0}^{m-1}
\sum_{x_1\in\mathcal{X}}\sum_{x_2\in\mathcal{X}}e^{2\pi i
h\frac{x_1-x_2}{m}}=m|\mathcal{X}|.
$$
We also note that for any $h,1\le h\le m/2,$ and any integers $L$
and $M\ge 1$ one has
$$
\left|\sum_{u=L+1}^{L+M}e^{2\pi
ih\frac{u}{m}}\right|=\frac{|\sin(\pi M/m)|}{|\sin(\pi h/m)|}\le
\frac{1}{|\sin(\pi h/m)|}\le \frac{m}{2h}.
$$

\section{A quantified version of the Weyl criterion}

Denote
$$
D(N,x_n)=\sup_{0\le \alpha<\beta\le 1}\left|\frac{F(N, x_n;
\alpha,\beta)}{N}-(\beta-\alpha)\right|.
$$
We describe our method in proving the following statement.

\begin{theorem}
\label{thm:LeVequetype} For any fixed real numbers $a$ and $b$ with
$a\ge 2b, 0\le b<2,$ the estimate
\begin{equation*}
D(N, x_n)\ll
\left(\sum_{h=1}^{\infty}h^{-\frac{2+a-2b}{2-b}}\left(\frac{1}{N}\left|\sum_{n=1}^Ne^{2\pi
ihx_n}\right|\right)^{\frac{2}{2-b}}\right)^{\frac{2-b}{2+a-b}}
\end{equation*}
holds, where the implied constant may depend only on $a$ and $b.$
\end{theorem}

In particular, taking $b=1$ one has for any fixed $a\ge 2$
$$
D(N,x_n)\ll
\left(\sum_{h=1}^{\infty}h^{-a}\left(\frac{1}{N}\left|\sum_{n=1}^Ne^{2\pi
ihx_n}\right|\right)^{2}\right)^{\frac{1}{a+1}}.
$$
If we take in the latter estimate $a=2,$  we obtain (apart from the
constant factor) LeVeque's inequality~\cite[p.9]{Monty}.

Taking $a=2b=4(1-\frac{1}{c}),$ one obtains for any fixed $c>1$
$$
D(N,x_n)\ll
\left(\sum_{h=1}^{\infty}h^{-c}\left(\frac{1}{N}\left|\sum_{n=1}^Ne^{2\pi
ihx_n}\right|\right)^{c}\right)^{\frac{1}{2c-1}}.
$$

Taking $a=2, b=2(1-\frac{1}{c}),$ one obtains for any fixed $c>1$
$$
D(N,x_n)\ll
\left(\sum_{h=1}^{\infty}h^{-2}\left(\frac{1}{N}\left|\sum_{n=1}^Ne^{2\pi
ihx_n}\right|\right)^{c}\right)^{\frac{1}{c+1}}.
$$

\proof It is easy to see that if we prove
\begin{equation}
\label{eqn:restrictedWeyl}
\frac{F(N, x_n; \alpha,\beta)}{N}-(\beta-\alpha)\ll
\left(\sum_{h=1}^{\infty}h^{-\frac{2+a-2b}{2-b}}\left(\frac{1}{N}\left|\sum_{n=1}^Ne^{2\pi
ihx_n}\right|\right)^{\frac{2}{2-b}}\right)^{\frac{2-b}{2+a-b}}
\end{equation}
in the case $1/4\le \beta-\alpha\le 1/2,$ then we are done. Indeed,
if $1/2\le \beta-\alpha\le 1,$ then~$1/4\le (\beta-\alpha)/2\le
1/2.$ Therefore,~\eqref{eqn:restrictedWeyl} can be applied to the
intervals
$$
[\alpha, \alpha+\frac{\beta-\alpha}{2})\quad {\rm and}\quad
[\alpha+\frac{\beta-\alpha}{2},\beta).
$$
This yields the required estimate for any $\alpha,\beta$ with
$1/2\le \beta-\alpha\le 1.$

If $0<\beta-\alpha<1/4,$ then consider the sequence $\{x_n\}-\alpha$
and apply~\eqref{eqn:restrictedWeyl}  with this sequence instead of
$x_n$ to the interval $[\beta-\alpha, 1).$ Then it remains to note
that
$$
F(N, x_n, \alpha,\beta)=N-F(N, \{x_n\}-\alpha; \beta-\alpha, 1)
$$ which follows from the fact that for any given $n, 1\le n\le N,$
either $\alpha\le\{x_n\}<\beta$ or
$\beta-\alpha\le\{\{x_n\}-\alpha\}<1.$

\bigskip

We now proceed to prove~\eqref{eqn:restrictedWeyl} for
$\alpha,\beta$ with $1/4\le \beta-\alpha\le 1/2.$ We may suppose
that $0\le x_n<1.$

Let us first reduce the problem to the case when $x_n$ are rational
numbers. Since
$$
W(N,
x_n):=\left(\sum_{h=1}^{\infty}h^{-\frac{2+a-2b}{2-b}}\left(\frac{1}{N}\left|\sum_{n=1}^Ne^{2\pi
ihx_n}\right|\right)^{\frac{2}{2-b}}\right)^{\frac{2-b}{2+a-b}}>0
$$
and since for any $L>10$
$$
\sum_{h=L+1}^{\infty}
h^{-\frac{2+a-2b}{2-b}}\left(\frac{1}{N}\left|\sum_{n=1}^Ne^{2\pi
ihx_n}\right|\right)^{\frac{2}{2-b}}\le L^{-\frac{a-b}{2-b}},
$$
then there exists a number $\delta>0$ such that for any sequence
$x'_n$ with the condition $|x'_j-x_j|\le \delta, j=1,\ldots, N,$ we
have
$$
|W(N, x'_n)-W(N, x_n)|\le W(N, x_n)/2.
$$
Thus
\begin{equation}
\label{eqn:WviaW} W(N, x'_n)<2W(N, x_n).
\end{equation}

Next, if for some $n\le N,$ $x_n\in [\alpha,\beta),$ then clearly we
can choose $x'_n$ to be a rational number such that
$$
x_n\le x'_n\le x_n+\delta,\qquad x_n'\in [\alpha, \beta).
$$
Besides, if $x_n\not\in [\alpha,\beta)$ then we can choose $x'_n$ to
be a rational number such that
$$
x_n\le x'_n < \min\{1, x_n+\delta\},\qquad x_n'\not\in [\alpha,
\beta).
$$
Hence, since any interval of positive length contains a rational
number, then we derive that there exists a sequence of rational
numbers $x'_n$ satisfying~\eqref{eqn:WviaW} and such that
$$
F(N, x_n; \alpha,\beta)=F(N, x'_n; \alpha,\beta).
$$
Thus, denoting $x'_n=s_n/m,$ where $s_n$ and $m>10$ are integers, we
conclude that it is indeed sufficient to prove the bound
$$
\frac{F(N, s_n/m; \alpha,\beta)}{N}-(\beta-\alpha)\ll
\left(\sum_{h=1}^{\infty}h^{-\frac{2+a-2b}{2-b}}\left(\frac{1}{N}\left|\sum_{n=1}^Ne^{2\pi
ihs_n/m}\right|\right)^{\frac{2}{2-b}}\right)^{\frac{2-b}{2+a-b}}.
$$
We can choose $m$ to be as large as we wish, just by substituting
$s_n/m$ by $ks_n/(km).$ In particular, we may assume that
$$
m^{1/2}W(N,s_n/m)>10,\quad m>(a+1)^2.
$$

Now observe that $F(N, s_n/m; \alpha,\beta)$ is equal to the number
of solutions of the congruence
$$
s_n\equiv y \pmod m,\quad n\le N,\quad \alpha m\le y<\beta m.
$$
Set
$$
R=\sum_{h=1}^{\infty}h^{-\frac{2+a-2b}{2-b}}\left(\frac{1}{N}\left|\sum_{n=1}^Ne^{2\pi
ihs_n/m}\right|\right)^{\frac{2}{2-b}}.
$$
If $R^{(2-b)/(2+a-b)}\ge 1/10,$ then the required estimate becomes
trivial. For this reason we suppose that $R^{(2-b)/(2+a-b)}<1/10.$
Take $k=[a]+1$ and define $T=[mR^{(2-b)/(2+a-b)}/k].$ Then
$$
kT <m/10<(\beta-\alpha)m/2,\quad (\beta-\alpha)m+kT<m, \quad T\ge
[10m^{1/2}/k]\ge 10.
$$
Let $J_1$ be the number of solutions of the congruence
$$
s_n\equiv y-y_1-\ldots-y_k\pmod m,
$$
where the variables are subject to the restriction
$$
n\le N,\quad \alpha m\le y< \beta m +kT, \quad 1\le
y_1,\ldots,y_k\le T.
$$
Here the length of the interval for $y$ is less than $(\beta-\alpha)
m+kT<m.$

Next, let $J_2$ be the number of solutions to the congruence
$$
s_n\equiv y+y_1+\ldots+y_k\pmod m,
$$
where the variables are subject to the restriction
$$
n\le N,\quad \alpha m\le y< \beta m -kT, \quad 1\le
y_1,\ldots,y_k\le T.
$$
Here, according to the choice of parameters we have $\alpha m<\beta
m-kT.$

Obviously
\begin{equation}
\label{eqn:twopolicenew} \frac{1}{T^k}J_2\le F(N, s_n/m;
\alpha,\beta)\le \frac{1}{T^k}J_1.
\end{equation}
Application of trigonometric sums yields
$$
\frac{J_1}{T^k}=\frac{1}{mT^k}\sum_{h=0}^{m-1}\sum_{n=1}^{N}\sum_{\alpha
m\le y< \beta m+kT}\sum_{y_1=1}^T\ldots\sum_{y_k=1}^{T}e^{2\pi i
h\frac{s_n-y+y_1+\ldots+y_k}{m}}.
$$
Picking up the term corresponding to $h=0$ and observing that for
$y$ there are $(\beta-\alpha)m+kT+\theta$ possible values, where
$|\theta|\le 1,$ we obtain
\begin{equation}
\label{eqn:asforJ1new}
\left|\frac{J_1}{T^k}-(\beta-\alpha)N\right|\le
\frac{2kTN}{m}+\frac{2}{mT^k}\sum_{1\le h\le
m/2}|S_1(h)||S_2(h)||S_3(h)|^k,
\end{equation}
where
$$
S_1(h)=\sum_{n=1}^{N}e^{2\pi i h\frac{s_n}{m}}, \quad
S_2(h)=\sum_{\alpha m\le y< \beta m+kT}e^{2\pi i h\frac{y}{m}},
$$
$$
S_3(h)=\sum_{y_1=1}^{T}e^{2\pi i h\frac{y_1}{m}}.
$$
Now we use the bound
$$
S_2(h)\ll m/h
$$
and also
$$
|S_3(h)|^{k}\le T^{k-a/2}|S_3(h)|^{a/2}\le
T^{k-a/2}\left(\frac{m}{h}\right)^{a/2-b}|S_3(h)|^{b}.
$$
Here we have used that $a\ge 2b.$ Incorporating this
into~\eqref{eqn:asforJ1new}, we obtain
$$
\left|\frac{J_1}{T^k}-(\beta-\alpha)N\right|\ll
\frac{TN}{m}+\frac{m^{a/2-b}}{T^{a/2}}\sum_{1\le h\le
m/2}h^{-1-a/2+b}|S_1(h)||S_3(h)|^b.
$$
Next, by Holder's inequality,
\begin{eqnarray*}
&&\sum_{1\le h\le m/2}h^{-1-a/2+b}|S_1(h)||S_3(h)|^b\le\\
&&\left(\sum_{h=1}^{\infty}h^{-\frac{2+a-2b}{2-b}}|S_1(h)|^{\frac{2}{2-b}}\right)^{(2-b)/2}
\left(\sum_{h=0}^{m-1}|S_3(h)|^2\right)^{b/2}=
\\
&&NR^{(2-b)/2}(mT)^{b/2}.
\end{eqnarray*}
Therefore,
$$
\left|\frac{J_1}{T^k}-(\beta-\alpha)N\right|\ll
\frac{TN}{m}+N\left(\frac{m}{T}\right)^{(a-b)/2}R^{(2-b)/2}.
$$
Recalling the choice of $T,$ we obtain
$$
\left|\frac{J_1}{T^k}-(\beta-\alpha)N\right|\ll NR^{(2-b)/(2-b+a)}.
$$

Analogously
$$
\left|\frac{J_2}{T^k}-(\beta-\alpha)N\right|\le NR^{(2-b)/(2-b+a)}.
$$
Therefore, from~\eqref{eqn:twopolicenew} we conclude that
$$
\left|F(N, s_n/m; \alpha,\beta)-(\beta-\alpha)N\right|\ll
NR^{(2-b)/(2-b+a)}.
$$
Theorem~\ref{thm:LeVequetype} is proved.

\section{Remarks}

Using the same argument one can deduce that if $0<\varepsilon\le 1,$
$\beta-\alpha\ge \frac{2\Delta}{\varepsilon}$ and if  the estimate
$$
\left|\sum_{n=1}^{N}e^{2\pi i hx_n}\right|\le \Delta N
$$
holds for any integer $h$ with $1\le h\le \Delta^{-1-\varepsilon},$
then
$$
F(N, x_n; \alpha,\beta)=(\beta-\alpha)N+O(\Delta
N\log\frac{\beta-\alpha}{\Delta}),
$$
where the implied constant in the $O-$symbol depends only on
$\varepsilon.$ This result does not follow from the Erd\H
os-Tur{\'a}n inequality, but it can be derived
from~\eqref{eqn:Monty}.

If one would like to have under hands only the proof of Weyl's
criterion, without its quantified version, then the argument given
in the previous section can be simplified even more. That is,
suppose that $0<\varepsilon<10^{-3}.$ We require the following
condition:
\begin{itemize}
\item[(i)] the inequality
$$
\left|\sum_{n=1}^Ne^{2\pi ihx_n}\right|\le \varepsilon^3 N
$$
holds for any integer $h, 1\le h\le \varepsilon^{-3}.$
\end{itemize}
Then we establish the following form of the Weyl criterion: {\it
under the condition {\rm(i)},
$$
\left|\frac{F(N, x_n; \alpha,\beta)}{N}-(\beta-\alpha)\right|\le
\varepsilon.
$$}
It is sufficient to show that
$$
\left|\frac{F(N, x_n;
\alpha,\beta)}{N}-(\beta-\alpha)\right|\le \varepsilon/2
$$
in the case $1/4\le \beta-\alpha\le 1/2.$ Then by continuity
argument the problem is reduced to the case with rational numbers,
that is for some integers $s_n$ and $m>100\varepsilon^{-1},$ we have
$$
\left|\frac{F(N, x_n;
\alpha,\beta)}{N}-(\beta-\alpha)\right|=\left|\frac{F(N, s_n/m;
\alpha,\beta)}{N}-(\beta-\alpha)\right|
$$
and
$$
\left|\sum_{n=1}^Ne^{2\pi ih\frac{s_n}{m}}\right|\le 2\varepsilon^3
N
$$
for any integer $h, 1\le h\le \varepsilon^{-3}.$

Now $F(N, s_n/m; \alpha,\beta)$ is equal to the number of solutions
of the congruence
$$
s_n\equiv y \pmod m,\quad n\le N,\quad \alpha m\le y<\beta m.
$$
Denote $T=[\varepsilon m/10]$ and set $J_1$ to be the number of
solutions of the congruence
$$
s_n\equiv y-y_1\pmod m, \quad n\le N,\quad  \alpha m\le y< \beta m
+T, \quad 1\le y_1\le T.
$$
Since $\beta-\alpha\le 1/2,$ then the length of the interval for $y$
is less than $m.$

Next, let $J_2$ be the number of solutions to the congruence
$$
s_n\equiv y+y_1\pmod m, \quad n\le N,\quad  \alpha m\le y< \beta m
-T, \quad 1\le y_1\le T.
$$
Since $\beta-\alpha\ge 1/4,$ then $\alpha m<\beta m-T.$

Obviously,
\begin{equation}
\label{eqn:twopolice} \frac{J_2}{T}\le F(N, s_n/m; \alpha,\beta)\le
\frac{J_1}{T}.
\end{equation}
For $J_1/T$ we have
$$
\frac{J_1}{T}=\frac{1}{mT}\sum_{h=0}^{m-1}\sum_{n=1}^{N}\sum_{\alpha
m\le y< \beta m+T}\sum_{y_1=1}^{T}e^{2\pi i h\frac{s_n-y+y_1}{m}}.
$$
Picking up the term corresponding to $h=0$ and observing that for
$y$ there are $(\beta-\alpha)m+T+\theta$ possible values, where
$|\theta|\le 1,$ we obtain
$$
\left|\frac{J_1}{T}-(\beta-\alpha)N\right|\le
\frac{2TN}{m}+\frac{2m}{T}\sum_{1\le h\le
m/2}h^{-2}\left|\sum_{n=1}^{N}e^{2\pi i h\frac{s_n}{m}}\right|.
$$
The sum over $h$ on the left hand side is
$$
\le 2\varepsilon^3N\sum_{1\le h\le
\varepsilon^{-3}}h^{-2}+N\sum_{h>\varepsilon^{-3}}h^{-2}\le
5\varepsilon^{3}N.
$$
Hence, recalling that $T=[\varepsilon m/10]$ and
$\varepsilon<10^{-3},$ we deduce
$$
\left|\frac{J_1}{T}-(\beta-\alpha)N\right|\le\varepsilon N/2.
$$

Analogously
$$
\left|\frac{J_2}{T}-(\beta-\alpha)N\right|\le\varepsilon N/2.
$$
Therefore, from~\eqref{eqn:twopolice} we conclude that
$$
\left|F(N, s_n/m; \alpha,\beta)-(\beta-\alpha)N\right|\le
\varepsilon N/2 .
$$

\bigskip

{\bf Acknowledgements.} This work was supported by Project
PAPIIT-IN105605 from the UNAM.

\end{document}